\documentclass[12pt]{article}

\usepackage{amsmath,amsthm,amsfonts,amssymb,epsf,epsfig}
\evensidemargin
.25in \textheight 8.5in \topmargin 0in 
\title{A Contamination Carrying Criterion for Branched Surfaces}
\author{ Ulrich Oertel \thanks{Research supported by  Max Planck Institute, Bonn; the National Science Foundationthe
National Science Foundation, grant DMS 9803293}
\hbox{\ \ }\&\hbox{\ }  
Jacek \'Swi\c atkowski \thanks{Research supported by  Max Planck Institute, Bonn; the National Science Foundation;  and
the Polish Scientific Research Committee (KBN), grant 5 P02A 035 20}}
\date{December, 2001}

\newtheorem{thm}{Theorem}[section] 
\newtheorem{lemma}[thm]{Lemma}

\theoremstyle{definition}
\newtheorem{defn}[thm]{Definition}

\newtheorem{example}[thm]{Example}

\theoremstyle{remark}
\begin{document}
\maketitle
%\tableofcontents

\def\reals {\hbox {\rm {R \kern -2.8ex I}\kern 1.15ex}}
\def\proj{P}
\def\hyp {\hbox {\rm {H \kern -2.8ex I}\kern 1.15ex}}
\def\integers {\hbox {\rm { Z \kern -2.8ex Z}\kern 1.15ex}}
\def\naturals {\hbox {\rm {N \kern -2.8ex I}\kern 1.15ex}}
\def\intr{\text{int}}
\def\inter{\ \raise4pt\hbox{$^\circ$}\kern -1.6ex}
\def\Cal{\cal}
\def\from{:}
\def\inverse{^{-1}}
\def\Max{\text{Max}}
\def\Min{\text{Min}}
\def\embed{\hookrightarrow}
\def\Genus{\text{Genus}}
\def\Z{Z}
\def\X{X}
\def\roster{\begin{enumerate}}
\def\endroster{\end{enumerate}}
\def\definition{\begin{defn}}
\def\enddefinition{\end{defn}}
\def\subhead{\subsection\{}
\def\theorem{thm}
\def\endsubhead{\}}
\def\head{\section\{}
\def\endhead{\}}
\def\example{\begin{ex}}
\def\endexample{\end{ex}}
\def\ves{\vs}
\def\mZ{{\mathbb Z}}
\def\M{M(\Phi)}
\def\bdry{\partial}
\def\hop{\vskip 0.15in}

\centerline{Rutgers University, Newark;  Wroc\l aw University, Wroc\l aw }

\abstract  A contamination in a 3-manifold is an object interpolating between the contact structure and the lamination. 
Contaminations seem to provide a link between 3-dimensional contact geometry and the classical topology of 3-manifolds, as
described in a separate paper \cite{UOJS:Contaminations}. 
In this paper we deal with contaminations carried by branched surfaces, giving a sufficient
condition for a branched surface to carry a pure contamination.
\endabstract

\section{Introduction }\label{Introduction}

\begin{figure}[ht]
\centering
\scalebox{1.0}{\includegraphics{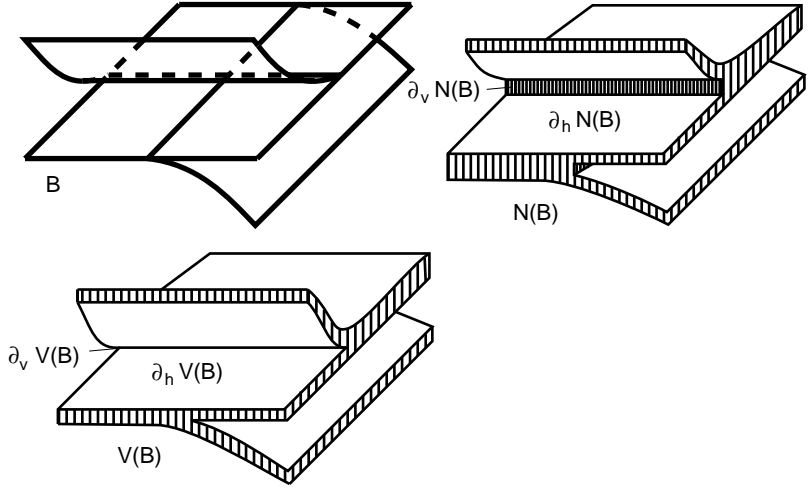}}
\caption{\small Local model for branched surface and fibered neighborhoods.}
\label{CarryNeighborhood}
\end{figure}

Let $M$ be an oriented 3-manifold and let $B\embed M$ be a closed branched
surface embedded in
$M$.  A contamination carried by B (which will be defined precisely later) is a plane field defined on a certain kind
of regular neighborhood of $B$.  In fact, we shall use two different
kinds of neighborhoods of branched surfaces.  Let $N(B)$ denote a fibered
neighborhood of the branched surface $B$.  It is foliated by interval
``fibers" and its boundary is decomposed into two parts: the vertical
boundary
$\bdry_vN(B)$ and the horizontal boundary $\bdry_hN(B)$, see Figure \ref{CarryNeighborhood}. 
There is a projection 
$\pi\from N(B)\to B$ which collapses the interval fibers foliating $N(B)$. 
The branched surfaces used in this paper have {\it generic branch locus},
meaning that the branched surface is locally modelled on the branched
surfaces shown in Figure \ref{CarryNeighborhood}.  Definitions related to branched surfaces can
be found in \cite{RW:BranchedSurfaces}, \cite{WFUO:IncompressibleViaBranched}, \cite{UO:MeasuredLaminations},
\cite{DGUO:EssentialLaminations}, and \cite{LMUO:Nonnegative}.   If we collapse the interval fibers of
$\bdry_vN(B)$, then we obtain another type of neighborhood,
$V(B)$ as shown in Figure \ref{CarryNeighborhood}.  Corresponding to $\bdry_v(N(B))$, which is
a union of annuli, we have 
$\bdry_v(V(B))$, which is a union of curves.   Again we have a projection,
also denoted $\pi$, which projects $V(B)$ to $B$, $\pi\from V(B)\to B$, and
which projects $\bdry_vV(B)$ to the branch locus of $B$.  Cutting $\bdry V(B)$
on the curves of $\bdry_vV(B)$, we obtain the horizontal boundary, $\bdry_hV(B)$
which is smoothly mapped to $M$.

A {\it positive contamination} $\xi$ carried by $B$ is a smooth plane field
defined on
$V(B)$ which is everywhere transverse to fibers of $V(B)$, is a positive
confoliation in $\text{int}(V(B))$, and is tangent to 
$\bdry_hV(B)$.  A negative contamination is defined similarly.  A {\it pure
positive contamination} carried by $B$ is a positive contamination which is 
contact in $\intr (V(B))$.  

\begin{figure}[ht]
\centering
\scalebox{1.0}{\includegraphics{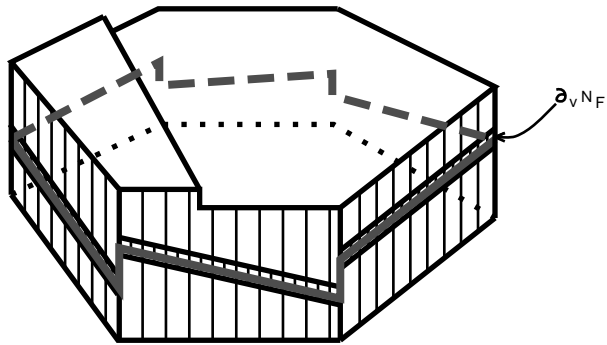}}
\caption{\small The pull-back $N_F$.}
\label{CarryN_F}
\end{figure}

For the following definition it will be convenient to choose a Riemannian metric
on a given generic branched surface $B$ such that branch curves intersect orthogonally.
A {\it positive twisted immersed surface of contact} for
$B$ is an immersion $f\from F\to N(B)$ mapping an oriented surface $F$ 
transverse to fibers of
$N(B)$ except possibly at $\bdry F$.  The map $f$ restricted to $\bdry F$
must map
$\partial F$ transverse to fibers in $\intr(\bdry_vN(B))$ except possibly on finitely
many closed disjoint intervals 
$C_1, C_2,\ldots, C_k$ in $\bdry F$, which are embedded by $f$ in fibers of
$N(B)$ corresponding to double points of the branch locus of
$B$.  We also require that $f(C_i)$ intersect $\intr (N(B))$, i.e. are not
contained in
$\bdry_vN(B)$.  We let
$C=\cup_iC_i$, and we say that each
$C_i$ is a {\it corner}.  The immersion must satisfy further conditions. 
We can pull back the fibered neighborhood structure to $F$ obtaining
$N_F$, a portion of a fibered neighborhood over $F$, see Figure \ref{CarryN_F},
homotopically equivalent to
$F$.  The orientation on $F$ and an orientation on $N_F$ determine an
orientation on fibers of
$N_F$. We require that $f$ embed each oriented $C_i$ to a fiber {\it
respecting orientations}:  The orientation induced on 
$\bdry F$ by the orientation on $F$ also gives an orientation for $C$, and $C$
must be mapped to $N_F$ such that the oriented arcs are mapped respecting
orientation to the oriented arcs of $N_F$.  For example, in Figure \ref{CarryN_F}, $N_F$ is shown for an immersed disk
of contact $F$, and the orientation on $F$ induces the upward orientation on fibers of $N_F$.  We see in the figure,
that the orientation on $C_i$'s induced by the orientation of $F$ are the same as the orientations induced by the
orientations of fibers of $N_F$.  There are more requirements.  The map $\pi\circ f$ descends to a map (also denoted
$f$) on the quotient space $F/C$, in which each arc
$C_i$ becomes a point which we call a {\it corner}. Pulling back the metric of $B$ to $F/\sim$ using the map
$f \from F/C\to B$, we require that at its corners $F/\sim$ have interior angles $\pi/2+2n\pi$ for some $n\ge 0$.  If
$C=\emptyset$, then we say $f$ is an immersed {\it surface of contact}, but not a
twisted immersed surface of contact; thus, for a twisted immersed surface of
contact we require
$C\ne
\emptyset$.

  A {\it negative twisted immersed 
surface of contact} is defined like the positive one, but now
$C_i$'s are mapped to interval fibers corresponding to double points of the
branch locus of
$B$ reversing orientations. 

\begin{figure}[ht]
\centering
\scalebox{1.0}{\includegraphics{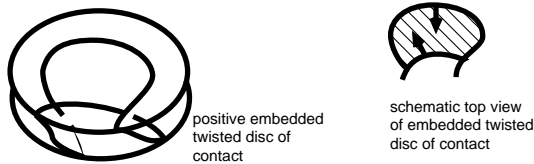}}
\caption{\small Twisted disc of contact sector.}
\label{CarryTwistedDisc}
\end{figure}

 If the map
$f\from F/C\to B$  is an embedding, it is easy to draw the
embedded twisted surface of contact as it appears in $B$.  We illustrate a
positive embedded twisted disk of contact in Figure \ref{CarryTwistedDisc}.  The figure includes
a schematic representation obtained by viewing the branch locus from
``above," where above is defined in terms of the transverse orientation of
the disk.  Of course, giving the disk the opposite orientation still gives
a positive twisted disk of contact of the same sign.  
In the schematic representation of a portion of branched surface, one also
needs to indicate the direction of branching.  If sectors $W,X,Y$ are
adjacent along an arc of branch locus
$\gamma$, and if
$W\cup Y$ and $X\cup Y$ are smooth, we say that branching along the arc
$\gamma\subset \bdry Y$ is {\it inward for} $Y$ and that branching along
the arc
$\gamma\subset\bdry X$ is {\it outward for $X$}.  We indicate the inward
direction with an arrow as shown in Figure \ref{CarryTwistedDisc}.

We should point out that Figure \ref{CarryTwistedDisc} is somewhat misleading, in that the
behavior at the boundary of a twisted immersed surface of contact can be
worse than illustrated.  Namely, at a corner $c_i$ of $F/C$, the image of
$F/C$ under $f$  may wrap around $f(c_i)$ (a double point
of the branch locus of $B$) more than one full turn, so that $f$
is not even locally an embedding in a neighborhood of $c_i\in F/C$.  This
is described in Section 2, see Figure \ref{CarryTiscModels}.

The name ``positive twisted immersed surface of contact" is too long, so we
reluctantly resort to the use of an acronym.  A {\it tisc} is a twisted
immersed surface of contact; a {\it positive tisc} is a positive immersed
twisted surface of contact; an {\it isc} is an immersed surface of contact. 

\begin{thm} \label{CarryTheorem}  Suppose $B\embed M$ is a branched surface. 
Suppose $B$ has no negative tiscs and no iscs.  Then
$B$ fully carries a positive pure contamination.
\end{thm}

{\it Open sectors} of a branched surface $B$ are the connected components
obtained after removing the branched locus from $B$; {\it sectors} are their completions
relative to a Riemannian metric on the branched surface.
A positive (negative) tisc for $B$ induces integer weights on
the sectors.  If $f\from F\to N(B)$ is a tisc,
the {\it weight} on a given sector is the number of components in the tisc of the
preimage of the open sector under the map $\pi\circ f$. The
weights assigned to all sectors give the {\it weight vector}. These weights must
satisfy certain equations and inequalities, called {\it positive (negative) tisc
equations and inequalities} which we shall describe in Section \ref{Weight}.  There are
similar  weight vectors for immersed surfaces of contact, which satisfy other relations
called {\it isc
equations and inequalities}. 

\begin{thm} \label{CarryTheoremWeights}   If a branched surface $B\embed M$ admits no
integer weight vectors satisfying the negative tisc equations and
inequalities, and it admits no integer weight vectors satisfying the isc equations and inequalities, then
$B$ fully carries a pure positive contamination.
\end{thm}

The proof of Theorem 1.1 occupies Sections 3-7. It is preceded
(in Section \ref{Weight})
by a discussion of weights and the proof of Theorem 1.2
as a corollary to Theorem 1.1.

\section{Weight Vectors and Proof of Theorem 1.2}\label{Weight}

In this section we discuss weight vectors induced by tiscs and iscs and prove
Theorem 1.2, assuming Theorem 1.1. 

\begin{figure}[ht]
\centering
\scalebox{1.0}{\includegraphics{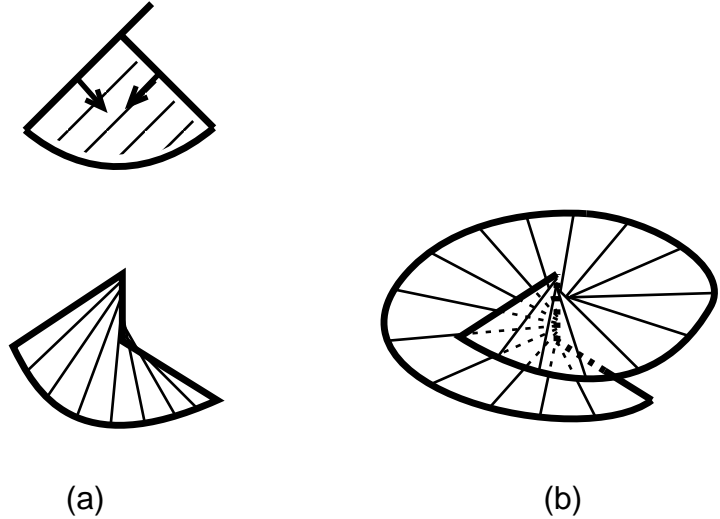}}
\caption{\small Corners of twisted discs of contact.}
\label{CarryTiscModels}
\end{figure}

As we mentioned in the introduction, a tisc can locally be more complex
near a corner than shown in Figure \ref{CarryTiscModels}a.  In general, a tisc can be
locally embedded in $N(B)$ as a truncated helix near  an interval fiber of
$N(B)$ which is the preimage of a double point of the branch locus, see
Figure \ref{CarryTiscModels}b.  The figure shows a helical surface rotating between 1 and 2
full turns; in general, any number of full turns can be added.  With a particular Riemannian 
metric on $B$ such that branch curves all intersect orthogonally, the angle of rotation after projecting to $B$ is
$\pi/2+2n\pi$ for some $n\ge 0$, and in the figure the angle is $5\pi/2$.  If a tisc has a corner with angle
$\pi/2+2n\pi$, $n>0$, on its boundary, it  cannot easily be represented using the schematic of Figure
\ref{CarryTiscModels}a.  However, we shall see that an arbitrary tisc can be replaced by one with ``convex" or $\pi/2$
corners only, and having the same weight vector.   Formally, if $f\from F\to N(B)$ is a tisc,
a {\it convex corner} is a corner of the surface $F/C$ such that $f\from F/C\to B$ 
is locally an embedding near this point.

\begin{figure}[ht]
\centering
\scalebox{1.0}{\includegraphics{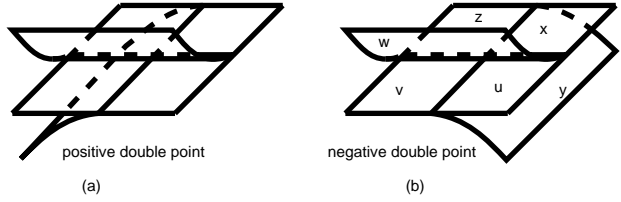}}
\caption{\small Weights on branched surface near double point of branch locus.}
\label{CarryWeights}
\end{figure}

In order to describe a practical method for detecting tiscs in a branched
surface $B\embed M$, we will describe the weight vectors induced by tiscs. 
If a branched surface has no weight vectors of this kind, we will conclude
that $B$ has no tiscs.  In a generic branched surface
$B\embed M$, there are two kinds of double points of the branch locus, as
shown in Figure \ref{CarryWeights}.  A positive (negative) double point is of the kind
that lies on the boundary of a positive (negative) tisc.  Looking at the
negative double point of Figure \ref{CarryWeights}, we can describe the integer weight
vectors corresponding to a negative tisc 
$f\from F\to N(B)$.  For each arc or closed curve of the branch locus 
with double points removed, the weights must satisfy inequalities like
$z\ge x+y$, called a {\it branch curve inequality}.  For this inequality, the curve
of the branch locus in question is an arc or closed curve common to the boundaries
of the sectors labelled $x,\ y,\ z$ respectively.  In the figure, we also have
branch curve inequalities
$z\ge w+v$,
$v\ge u+y$, and $x\ge w+u$.  Each inequality is strict if a
portion of $\bdry F$ is mapped by $\pi\circ f$ to the edge corresponding to
the inequality.  The boundary of a {\it positive} tisc cannot turn the
corner at a negative double point, so for a positive tisc, we would also
have the {\it corner equation}
$$z-(w+v)=x-(u+w) \text{ or } z-v=x-u.$$   For a negative tisc, there are
similar equations at positive double points. For a negative tisc at a
negative double point, we have a {\it corner inequality}
$$z-(w+v)\ge x-(u+w).$$ Finally, for a negative tisc there must exist a
double point where the above inequality is strict.  Of course, all weights must be $\ge 0$.  The equations and
inequalities satisfied by weight vectors induced by a negative (positive) tisc will be referred to as
the {\it negative (positive) tisc equations and inequalities}.   An integer weight vector satisfying these
equations and inequalities is called a {\it negative (positive) tisc
weight vector}. 

Similarly, an immersed surface of contact or {\it isc}, connected or not,
determines an isc weight vector (with at least one curve inequality strict and
all corner equations satisfied).  The equations and inequalites satisfied by weights induced by
an isc are called {isc
equations and inequalities}.  An integer weight vector satisfying these equations and inequalities is called an
{\it isc weight vector}.

\begin{lemma}\label{TiscEquivalentTiscVector}  A negative (positive) tisc $f\from F\to N(B)$ in a branched surface
$B\embed M$ induces a negative (positive) tisc weight vector.  Conversely given a negative
(positive) tisc weight vector on $B$ (satisfying at least one strict corner
inequality) there is a surface $F$ and an immersion $f\from F\to N(B)$ such that the restriction to at
least one component is a negative (positive) tisc, with possible closed and surface of contact
components.  The tisc can be constructed such that every corner is convex.

A surface of contact $f\from F\to N(B)$ induces a weight vector satisfying isc
equations and inequalities (including at least one strict branch curve
inequality).  Conversely, for any weight vector satisfying the isc equations and
inequalities there is a surface $F$ and a carrying map $f\from F\to N(B)$ such that the restriction to at
least one of the components is a surface of contact, with possible closed surface
components.
\end{lemma}

\begin{proof}
We have already proved everything except the statement that a weight vector
satisfying negative (positive) tisc equations and inequalities is induced by a
surface at least one of whose components is a tisc, with all non-tisc 
components being surfaces
of contact or closed surfaces, and with only convex corners.  This is proved by
showing that the appropriate numbers of copies of each sector in $N(B)$ can be
glued at their edges and corners to yield the required tiscs and iscs, allowing
self intersections of the surface.  

Similarly, a weight vector satisfying isc
equations and inequalities is induced by a surface at least one of whose
components is an isc, and possibly also including closed components.
\end{proof}

One might then ask whether a negative (positive) tisc weight vector
determines a tisc
$f\from F\to N(B)$ {\it embedded} in
$N(B)$, but with
$f\from F/C\to B$ possibly not an embedding.  (Recall that invariant weight
vectors, which satisfy all switch equations, uniquely determine measured
laminations carried by $B$, and uniquely determine surfaces carried by $B$
if all weights are non-negative integers.)  

\begin{lemma} \label{TiscVectorGivesEmbeddedTisc}  A negative (positive) tisc (isc) weight vector $w$ for a branched surface
$B\embed M^3$ determines an embedding
$f\from F\to N(B)$ whose restriction to at least one component is an embedded tisc
(isc).  The tisc (isc) is unique up to isotopy through carrying maps, and in general
such a tisc has non-convex corners.
\end{lemma}

\begin{proof}  According to Lemma \ref{TiscEquivalentTiscVector}, a negative (positive) tisc (isc) weight
vector determines and immersion $f\from F\to N(B)$, which is a tisc (isc) on at
least one connected component.  Putting the map $f$ in general position, but still
transverse to fibers of $N(B)$, we can then perform cut-and-paste on curves of
self-intersection to obtain an embedding.  This may replace convex corners by
non-convex ones (as in Figure \ref{CarryTiscModels}b).
\end{proof}

\begin{lemma} \label{TiscGivesConvexTisc}  For any negative (positive) tisc $f_1\from F_1\to N(B)$, there is
another negative (positive) tisc $f_2\from F_2\to N(B)$ with the same weight vector and only convex corners.
\end{lemma}

\begin{proof}  From a negative tisc we obtain a negative tisc weight vector.  From the negative tisc
weight vector, Lemma \ref{TiscEquivalentTiscVector} gives an immersed negative tisc with convex corners.
\end{proof}

\begin{proof} {(\it Theorem 1.2.)}
The theorem follows by combining Theorem 1.1 with Lemma 2.1.
\end{proof}

\section{Elementary Splitting Moves}\label{Moves} 

The proof of Theorem 1.1 occupies Sections 3-7. 
Its first essential part (described in Section 4) consists of performing a sequence 
of splittings of the branched surface $B$, see \cite{UO:MeasuredLaminations},
\cite{DGUO:EssentialLaminations}, and \cite{LMUO:Nonnegative}, for precise
definitions of branched surface splitting. In this section we describe
(and discuss some properties of) elementary
splittings of the kind we will use in the next section.

\begin{figure}[ht]
\centering
\scalebox{1.0}{\includegraphics{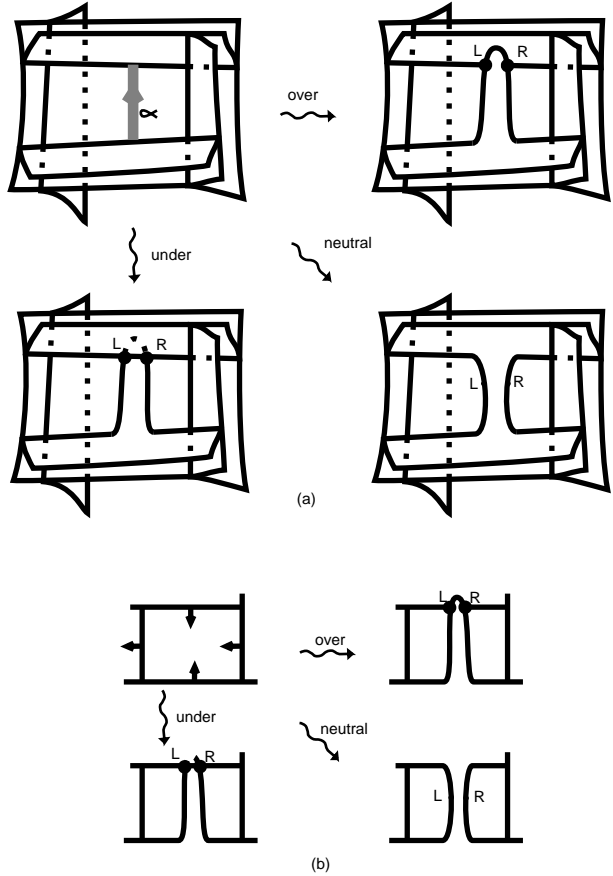}}
\caption{\small Choices for splitting.}
\label{CarryChoices}
\end{figure}

The typical splitting we use focuses on a particular
sector $Z$ of $B$ which has at least some inward branch locus on its
boundary, as in an embedded twisted disk of contact.  A {\it good directed
arc} properly embedded in $Z$ is a directed arc $\alpha$ with starting
point $P$ in $\bdry Z$ and ending point $Q$ in $\bdry Z$ which has the
property that the branch locus at $P$ and $Q$ is inward for $Z$.   
Starting at $P$, we split the branched surface as shown in Figure \ref{CarryChoices}.  When
we arrive at $Q$, we have a choice as to how we do the splitting at $Q$. 
The first part of the branch locus may pass ``over" or ``under" the other,
see Figure \ref{CarryChoices}. If we choose the {\it over move}, we call the branched
surface resulting from the move $B_o$; if we choose the {\it under move} we
call the resulting branched surface $B_u$.  A third possibility is for the
two segments of branch locus to meet as shown in Figure \ref{CarryChoices}, in a move that
we call the {\it neutral move}, resulting in a branched surface
$B_n$.  In this paper, we do not need to use the neutral move.  Figure \ref{CarryChoices}b shows the same
moves using the schematic representation of a sector of
$B$.  We shall contrive always to split
$B$ on good directed arcs passing through sectors of
$B$.

At the level of the fibered neighborhood $N(B)$ the splitting move can be
thought of as the removal of an $I$-bundle from $N(B)$.  If $B'$ is the
branched surface obtained after the splitting, then $N(B)=N(B')\cup J$,
where $J$ is an $I$-bundle over a surface.  In our situation $J=D\times I$,
where $D$ is a disk.  There is an arc $\beta\subset\bdry D$ such that
$J\cap \bdry_vN(B)=\beta\times I$.  Thus we are ``cutting" some of the
fibers of $N(B)$ where they meet a disk embedded in $N(B)$ transverse to
fibers.  

\begin{lemma}\label{SplittingLemma}  Suppose $B\embed M$ is a branched surface with
generic branch locus and without iscs or negative tiscs.  Suppose
$Z$ is a sector of $B$ and suppose $\alpha$ is a good directed arc in $Z$,
with
$\bdry\alpha\subset \bdry Z$, and with the beginning at $P$.  One of the
branched surfaces
$B_o$ or $B_u$, obtained by splitting along the arc $\alpha$ (using
the over or under move) is a branched surface without iscs or negative
tiscs.
\end{lemma}

\begin{proof} We consider moves in which one arc of the branch locus
passes over or under the other (Figure \ref{CarryChoices}); thus two new double points of
the branch locus are produced.  We label these $L$ (left) and $R$ (right),
whether we do the over-move or the under-move. If there is a new tisc, its
boundary must include at least one of $L$ or $R$ at least once, but possibly
more often.  For each of the two choices, over or under, for the move, a
new negative tisc might be produced, but the negativity shows that only one
of the two double points can be included on the boundary of the new
negative tisc.   After the over-move, a negative tisc can only include the
new double point $L$, while after the under-move a negative tisc can only
include the double point $R$, see Figure \ref{CarryInductionStep}.  (Similarly if a positive tisc
is produced after the under-move, it can only include $L$.)

\begin{figure}[ht]
\centering
\scalebox{1.0}{\includegraphics{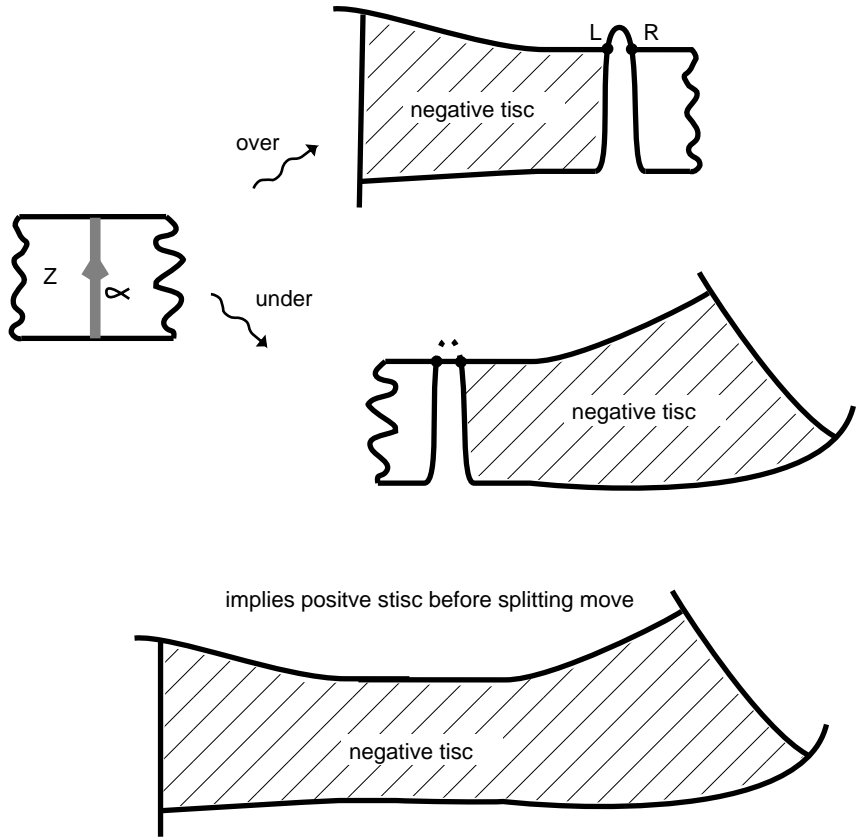}}
\caption{\small The induction step.}
\label{CarryInductionStep}
\end{figure}

If one of the moves produces no iscs and no negative non-trivial tiscs, this is the
move we choose to continue our splitting.  We must show that at least one
of the moves produces no iscs or negative tiscs.  Suppose the
over-move produces a negative tisc on the left, and the under-move produces
a negative tisc on the right.  Combining these using a boundary connected
sum operation yields a negative tisc before the move as shown in Figure \ref{CarryInductionStep},
which contradicts our induction hypothesis.  Of course the diagrams deal
only with one very special case, but the construction is roughly the same
in the general case.  Recall that by Lemma \ref{TiscEquivalentTiscVector}, we can restrict attention to tiscs with
convex corners.  However, it is possible that each of the two negative tiscs after the move has one
of the two double points $L$ or
$R$ represented on its boundary {\it more than once}.  For this reason, one
may have to use multiple copies of each tisc to obtain the tisc before the
move.  More precisely, if the negative tisc $f_1\from F_1\to N(B)$ on the
left produced after the over-move passes through the new double point $L$
on the left $\ell$ times, and if the negative tisc $f_2\from F_2\to N(B)$
on the right produced after the under-move passes through the double point $R$ on the right
$r$ times, then we produce the tisc  $f\from F\to N(B)$ by combining
$\text{lcm}(\ell,r)/\ell$ of the tiscs $f_1$ with $\text{lcm}(\ell,r)/r$
of the tiscs $f_2$.  (The fact that we may have to join many copies of the
two tiscs explains why ruling out immersed twisted {\it disks} of contact
in the assumptions of Theorem 1.1 appears not to be enough.) 

There is a slight error in the above argument.  If both the over and under-move produce a negative
tisc, these can be combined to form a new tisc before the move, but we must
ensure that there are corners (twisting) on 
$\bdry F$ which is obtained by joining the boundaries of the two tiscs
after the move.  This is not the case if for all possible new positive
tiscs produced by the over- or under-move, all corners on boundary
components intersecting a new double point of the branch locus are mapped
to the new double points ($L$ and $R$).  However, in this case it is easy
to see that there would have been an isc in the branched surface before 
modification.

It is also easy to see that the move cannot introduce an isc.

Note that in this proof we never use the neutral move.
\end{proof}

There is another move that one could use.  Suppose that in Figure \ref{CarryChoices}a the
directed arc $\alpha$ had its final endpoint on an arc of the branch locus where the
branching is outward with respect to the sector containing $\alpha$, i.e. with the
opposite sense of branching.  Then there is no choice for how to split on a neighborhood
of the arc $\alpha$.  In such a move, a portion of one arc of the branch locus
``passes" another arc with the same sense of branching.  It is best to avoid this kind of move, and we call the move a
{\it bad move}. The move is undesirable, because it can easily result in a branched surface containing a negative
tisc.

\section{Lamination with negative holonomy}\label{Lamination}

Our idea for proving Theorem 1.1 is first to construct an
auxiliary object $\Lambda$, a lamination in part of $V(B)$,
which in some sense encodes a sufficient amount of the ``twisting"
seen in contact plane fields, and which can then later be 
uniformly distributed over the interior of $V(B)$. The twisting
manifests itself in the lamination as a property, explained below, that $\Lambda$
has strictly negative holonomy. Conversion of $\Lambda$ into
a pure contamination in $V(B)$ is a technically involved but
fairly standard procedure, the details of which occupy later
sections.  We regard the construction of the lamination $\Lambda$, described in this
section, as the most delicate part of the proof of Theorem 1.1.

We will obtain $\Lambda$ as an inverse limit of a sequence of branched surfaces
obtained from $B$ by elementary splitting operations.  We
start by showing,  under assumptions of Theorem 1.1,
how to choose a sequence of elementary splittings appropriately.

We are given a branched surface $B\embed M$
with generic branch locus, and without negative tiscs.  We choose for $B$ a
structure as a 2-complex, where the branch locus is a subset of the
1-skeleton $X=B^{(1)}$ and the double points of branch locus are contained
in the 0-skeleton $B^{(0)}=X^{(0)}$.  We choose a reasonable metric for $B$
or $M$ such that if $\epsilon$ is sufficiently small then the 
$\epsilon$-neighborhood (or smaller) of $X$
gives a regular neighborhood of $X$.  Now choose a decreasing
sequence $\epsilon_n,\ n\ge 0$ of small numbers with
$\epsilon/2<\epsilon_n\le\epsilon$.

\begin{figure}[ht]
\centering
\scalebox{1.0}{\includegraphics{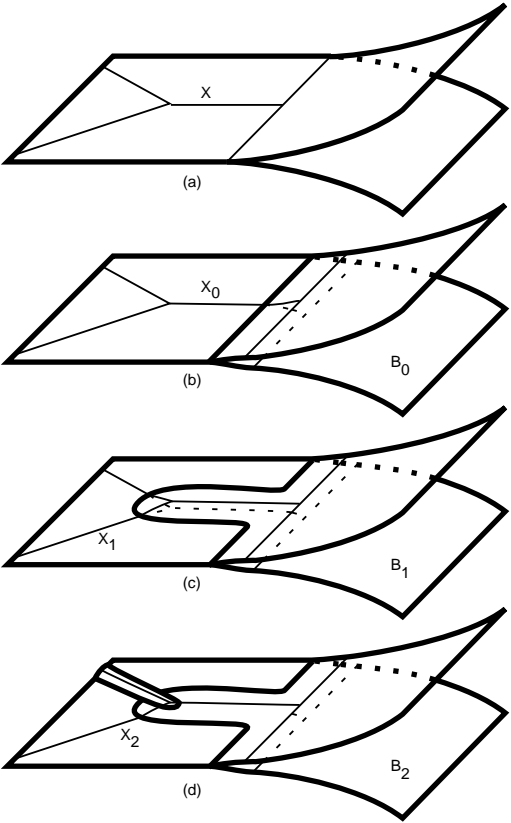}}
\caption{\small Splitting near the 1-complex.}
\label{CarrySplitOnComplex}
\end{figure}

When $n=0$, we begin with $B$ and split in an $\epsilon_0$ neighborhood of
the branch locus of $B$, see Figure \ref{CarrySplitOnComplex}ab.  This gives a branched surface
$B_0$ isomorphic to $B$, but now the 1-complex
$X$ pulls back, under the projection $\pi_0\from B_0\to B$ to a more
complex 1-complex 
$X_0\subset B_0$.  In particular, the pull-back of some edges of $X$ will
be train tracks.  We choose a directed 1-cell
$\alpha$ of
$X_0$ emanating from the branch locus of $B_0$ and split in an $\epsilon_1$
neighborhood of the arc as shown in Figure \ref{CarrySplitOnComplex}bc.  Note that $\alpha$
corresponds to a directed arc in $X$.  In fact, $\alpha$ projects to a
subarc of a 1-cell in $X$, with at most length $\epsilon$ truncated from
each end. It may happen that the final point of
$\alpha$ is also on the  branch locus of
$B_0$, so that $\alpha$ is a good directed arc for $B_0$.  In that case, we
must use Lemma \ref{SplittingLemma} to decide whether the newly formed branch locus goes
over or under the branch locus at the end of the arc.  The correct
choice ensures that the new branched surface
$B_1$ which we obtain will have no non-trivial negative tiscs.

We will
continue splitting along (pull-backs of) other arcs of
$X$, one at a time, to obtain a sequence $B_n$ of branched surfaces, each
projecting to $B$ via $\pi_n\from B_n\to B$.  At the
$n$-th step we split on a strip corresponding to an $\epsilon_n$
neighborhood of an arc contained in the pull-back train  track of a 1-cell
of $X$.  Note that later splittings are done on thinner strips, older
splittings are  done on fatter strips as shown in Figure \ref{CarrySplitOnComplex}. Figure \ref{CarrySplitOnComplex}cd
shows a splitting on a good arc, where Lemma \ref{SplittingLemma} is applied to make a
choice among possible splittings.  At the
$n$-th step, we choose to split
$B_n$ on an arc emanating from the branch locus near a 0-cell, and we do
this at {\it the oldest part} of the branch locus.  This means that we
begin the splitting at an available part of the branch locus {\it
outermost} in a neighborhood of a 0-cell of $X$.  This, 
and the fact that the sequence $\epsilon_n$ is decreasing, guarantees that
none of the new branch locus ``passes" the old branch locus to give a bad
move.  At the final end of the arc of splitting, however, we may need to
apply Lemma \ref{SplittingLemma} repeatedly to decide whether the splitting goes over or
under the opposing branch locus.  The correct choices ensure that
$B_n$ has no iscs or negative tiscs.  

\begin{figure}[ht]
\centering
\scalebox{1.0}{\includegraphics{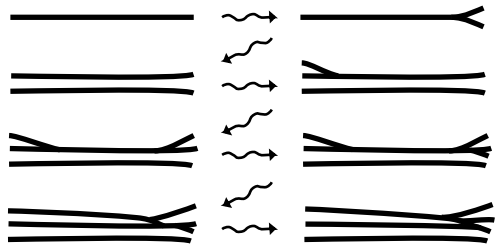}}
\caption{\small The effect of splitting.}
\label{CarrySplitSegment}
\end{figure}

We have projections $p_n\from B_n\to B$.  We illustrate a possible sequence
of events locally by choosing an edge $e$ in $X$, then examining
$p_n\inverse(e)\subset X_n$.  Figure \ref{CarrySplitSegment} shows such a sequence of splittings
of $e$ induced by the splittings of $B$.  The first three splittings are
consistent with the splittings induced on the central edge of
$X$ in Figure \ref{CarrySplitOnComplex}.

Denote by $N(X)$ the $\epsilon/2$-neighbourhood of $X$ in $B$. 
The infinite sequence of splittings as above defines a lamination $\Lambda$ 
in the preimage $\pi^{-1}(N(X))\subset V(B)$ as an inverse limit of the sequence of
projections $q_n\from B_j\to B_{j-1}$
$$\cdots  B_4\rightarrow_{q_4} B_3\rightarrow_{q_3} B_2\rightarrow_{q_2} B_1\rightarrow_{q_1} B_0\rightarrow_{q_0}
B.$$ 
The inverse limit can be realized as a space embedded in $M$, where it is the intersection of appropriately
chosen nested neighborhoods $N(B_j)\embed M$.  See \cite{LMUO:Nonnegative} for more details on
constructing laminations as inverse limits.

We can describe the inverse limit in another way.  For each splitting in the sequence (including the first splitting
from
$B$ to $B_0$) consider a sheet of surface in $V(B)$ ``parallel'' to this
splitting, i.e. a sheet of surface transverse
to the vertical fibers in $V(B)$, with boundary equal to the union of
the old and new part of the branch locus corresponding to the splitting.
We may think of this sheet as the splitting surface.  
The union of all such surfaces corresponding to the splittings in the infinite sequence,
intersected with the preimage $\pi^{-1}(N(X))$ is a non-compact splitting surface which we call $L$. The lamination
$\Lambda$ is obtained from $\pi\inverse(N(X))\subset V(B)$ by splitting on $L$.    Each
component of $L$ has boundary, with one arc of the boundary attached to the branch locus $\partial_v
V(B)$ and the remainder of the boundary in mapped by
$\pi$ to the boundary of $N(X)$ in $B$. This follows from the choices of splittings at the
oldest parts of branch locus, which guarantees that no part of branch
locus in the interior of $\pi^{-1}(N(X))$ remains unaffected by later
splittings.  Thus we can say that $L$ is carried by $N(X)$.

Crucial for our purposes is the following
property of the lamination $\Lambda$ or the splitting surface $L$.
Let $E$ be a 2-cell of $B$ and $\hat E$
the smaller 2-cell obtained by removing the part of $N(X)$ contained in $E$. 
Then the intersection with $L$ defines a nonempty immersed 1-dimensional
manifold in
the annulus $\pi\inverse(\bdry \hat E)=\bdry\hat E\times I$ in $V(B)$.  The intersection
with $\Lambda$ gives a 1-dimensional lamination. 

By the construction of $L$ and $\Lambda$, both $L$ and $\Lambda$ have
{\it strictly negative holonomy}.  This means that following any leaf
of $L$ or $\Lambda$ in the direction of the orientation induced on
$\bdry\hat E$ from the orientation of $B$, after a full turn around
$\hat E$, we end at a point in the same vertical fiber lying 
below the initial point using the canonical orientation
of the fibers.  (It is an unfortunate fact of life in 
contact geometry that positive contact structures induce negative holonomy.)  To prove this property 
observe that if it
were not true then we would have either a disk of contact (coming from
a leaf without holonomy) or a negative twisted disc
of contact (coming from a leaf
with positive holonomy) for some branched surface $B_n$ in our sequence.

We will refer to the above property by saying that the splitting surface $L$ or the lamination $\Lambda$
itself has {\it strictly negative holonomy}.

We end this section by noting that the arguments can be extended, without too much difficulty, to the case of a
noncompact branched surface $B$.

\section{From Lamination to Foliation}\label{Foliation}

In this section we start converting the splitting surface $L$ or the lamination $\Lambda$
with strictly negative holonomy into a pure
contamination in $V(B)$. As a first step we convert $L$ or $\Lambda$
to a smooth foliation $\cal F$ in $\pi^{-1}(N(X))$.

Without loss of generality we may assume $L$ has the 
following two further properties:
\item{(1)} $L$ is continuous, i.e. the tangent plane $T_xL$
depends continuously
on the point $x\in L$ (in the topology on $L$ induced
from $V(B)$);
\item{(2)} for each vertical fiber $I=[a,b]$ in $\pi^{-1}(N(X))$
we have $\sup L\cap I=b$ and $\inf L\cap I=a$.

\hop

\noindent
These two properties make it possible to extend $L$, regarded as a continuous
plane field defined only at some points of $\pi^{-1}(N(X))$, to a continuous integrable 
plane field defined throughout $\pi^{-1}(N(X))$. This gives a 
foliation $\cal F$ in $\pi^{-1}(N(X))$, transverse to
the fibers, having strictly negative holonomy at each sidewall
$\partial \hat E\times I$ (except at $\partial E\times\partial I$, i.e.
at the top and bottom). 
The foliation $\cal F$ typically contains singular leaves which are non-compact branched surfaces
with no double points in the branch locus, and which intersect $\bdry_hV(B)$.

A continuous foliation $\cal F$ with singular branched leaves with the required properties can also be constructed from
the lamination
$\Lambda$ roughly as follows.  One replaces boundary leaves of $\Lambda$ by product families of leaves, suitably
tapered, then one collapses gaps in $V(B)-\Lambda$.

Our next goal is to modify $\cal F$ so that 
it becomes smooth, and still has strictly negative holonomy.

Using a fiber-preserving homeomorphism of $V(B)$ we can make 
$\cal F$ smooth near the vertex fibers (i.e. vertical fibers
corresponding to vertices of $X$).
We can then choose smooth cylindrical charts 
$$\left\{(r,\theta,z):r\le1,\theta\in S^1, z\in I\right\}$$
in cylinders $E\times I$ with the following properties:

\item{(1)} lines $(r,\theta)=\text{constant}$ correspond to vertical fibers of $V(B)$;
\item{(2)} for some $\delta>0$, $\pi\inverse(N_{\epsilon/2}X)\cap E=\{(r,\theta,z):r\ge
1-\delta\}$;
\item{(3)} $\cal F$ is given in coordinates $(r,\theta,z)$ as the
kernel of a 1-form $dz+f(z,\theta)d\theta$, with $f$ independent of $r$,
$f\le 0$, and $f=0$ only near the arguments $\theta$ corresponding
to the vertex fibers in $\partial E\times I$ (where $\cal F$ is smooth).

Consider a 1-cell $A$ of $X$ which is the intersection of two 2-cells
$E_1$ and $E_2$ whose union $E_1\cup E_2$
is smooth at $A$. We focus on the part
of the foliation $\cal F$ restricted to $\pi^{-1}(A)$ (or more
precisely, to the intersection of the cylinders in $V(B)$ corresponding
to $E_1$ and $E_2$ in the preimage $\pi^{-1}(A)$).
In both charts $(\theta_1,z_1)$ and $(\theta_2,z_2)$
(in $\partial E_1\times I_1$ and $\partial E_2\times I_2$ 
respectively), this restriction is described by the continuous
functions $f_1$ and $f_2$ as in condition (3) above.
We will modify this part of $\cal F$ so that the functions $f_1$
and $f_2$ become smooth and still satisfy the assertions of condition
(3).

Let $f_1'\ge f_1$ be a nonpositive smooth function with the same
support as $f_1$ (such a function clearly exists). Modify $\cal F$
inside $E_1\times I_1$ and $E_2\times I_2$ (only close to $\pi^{-1}(A)$)
by making it equal to the kernel of the 1-form 
$dz_1+f_1'(\theta_1,z_1)d\theta_1$
and the 1-form $dz_2+f_2'(\theta_2,z_2)d\theta_2$, 
where $f_2'$ is the smooth function induced by $f_1'$ and the compatibility
of modifications at the intersection of the two charts.

The modified foliation still has strictly negative holonomy at 
$\partial E_1\times I_1$. The holonomy function is however smaller in absolute value
(or at least not bigger) than before modification. Due to the fact that the
orientations on $A$ induced by coordinates $\theta_1$ and $\theta_2$ are opposite,
we have $f_2'\le f_2$, which clearly implies that the holonomy
at $\partial E_2\times I_2$ of the modified foliation is also
strictly negative. Since the modification did not affect other
cylinders, the holonomy of the modified foliation is strictly
negative for all of them. Performing this sort of modification
at all edges $A$ in $X$ we get a smooth foliation $\cal F$ as
required.

\section{Special Charts}\label{Charts}

To prove Theorem 1.1 we will further modify the foliation $\cal F$ obtained
in the previous section. For this modification, which is presented in Section 7,
we need some coordinate charts well suited to the contamination structure.
This section contains description of such charts, which are similar to ones used
in \cite{UOJS:Contaminations}.  

Let $\xi$ denote a contamination
carried by $B\embed M$.

\hop

\begin{definition} Let $(r,\theta,z)$ be the cylindrical coordinates
in $R^3$, and let ${\cal C}=\{(r,\theta,z)\in R^3:r<R,-1\le z\le 1\}$. A {\it $\cal
C$-chart} or {\it cylinder chart} for a contamination $\xi$ carried by
$B\embed M$ is a smooth embedding $\psi:{\cal C}\to V(B)$, such that:
\item{(1)} the images in $M$ of all curves $\{r=const,\theta=const\}$
in $\cal C$ are contained in fibers of $V(B)$;
\item{(2)} the images in $M$ of all curves $\{\theta=const,z=const\}$
in $\cal C$ are everywhere tangent to $\xi$;
\item{(3)} the images in $M$ of the disks $\{z=\pm1\}$ in $\cal C$ are contained in the horizontal
boundary of $V(B)$.

\end{definition}

Observe that in view of (3), the image curves in condition (1) coincide with the fibers
of $V(B)$.

\hop

We mention without
proof the following easy fact.

\begin{lemma} \label{CChartLemma} Let $\xi$ be a positive contamination carried by $B\embed M$.  For any
open disk $D$ contained in the interior of a sector of $B$ (i.e. not intersecting the
branch locus) and for any smooth radial coordinates $(r,\theta)$ in $D$
there exists a $\cal C$-chart $\psi:{\cal C}\to V(B)$ such that
$\pi(\psi({\cal C}))=D$ and $\pi\circ\psi(r,\theta,z)=(r,\theta)$.

\end{lemma}

\begin{definition}\label{BChartDef} Let $(x,y,z)$ be the cartesian coordinates
in $R^3$, and let ${\cal B}=\{(x,y,z)\in R^3:x,y\in (-1,1),\ z\in[-1,1]\}$.
A {\it $\cal B$-chart or box chart} for a contamination $\xi$ carried by $B\embed M$
is a smooth embedding $\psi:{\cal B}\to V(B)$, such that
\item{(1)} the images in $M$ of all curves $\{x=const,y=const\}$
in $\cal B$ are contained in fibers of $V(B)$;
\item{(2)} the images in $M$ of all curves $\{x=const,z=const\}$ in $\cal B$ are everywhere tangent to $\xi$;
\item{(3)} the images in $M$ of the subsets $\{z=\pm1\}\subset{\cal B}$ are contained in
the horizontal boundary.

In some situations it is better to use a chart with $z\in (-1,1)$ and with
condition (3) above omitted. We will call a chart of this sort
an {\it open $\cal B$-chart}.
\end{definition}

The freedom for constructing $\cal B$-charts for a contamination $\xi$
is similar to that for $\cal C$-charts, see Lemma \ref{CChartLemma}. We omit the details.

\hop

The advantage of working with $\cal C$-charts and $\cal B$-charts
for a contamination $\xi$ carried by $B$ is that in the coordinates provided
by these charts $\xi$ can be expressed in a unique way
as the kernel of a 1-form $dz+f(r,\theta,z)d\theta$,
or $dz+f(x,y,z)dx$, 
for some smooth function $f$. We will call the function $f$ as above a {\it slope function}. 
It is possible to modify a contamination $\xi$ inside
a chart neighbourhood by modifying the appropriate slope function.
The following lemmas show how to express the property of being a confoliation 
(or a contact structure) in terms of a slope function.

\begin{lemma}\label{BoxChartLemma}  Let $\omega=dz+f(x,y,z)dx$ be a 1-form
and $\xi=\ker(\omega)$ be the induced plane field in $\cal B$.
Then $\xi$ is a positive confoliation if and only if $\partial f/\partial y(x,y,z)\le 0$
for each $(x,y,z)\in{\cal B}$. Moreover, $\xi$ is contact at a point $(x,y,z)\in{\cal B}$ if and only if $\partial
f/\partial y(x,y,z)<0$.
\end{lemma}

\begin{lemma}\label{CylinderChartLemma} Let $\omega=dz+f(r,\theta,z)d\theta$ be a 1-form
and $\xi=\ker(\omega)$ be the induced plane field in $\cal C$.
\item{(a)} The form $\omega$ is well defined and smooth in $\cal C$ iff
$f(r,\theta,z)=r^2\cdot h(r,\theta,z)$ for some smooth function
$h:{\cal C}\to R$.
\item{(b)} The plane field $\xi$ is a positive confoliation iff
${{\partial f}/{\partial r}}(r,\theta,z)\le 0$ for all points in $\cal C$ with $r>0$.
\item{(c)} The plane field $\xi$ is a positive contact structure
at a point $(r,\theta,z)\in{\cal C}$ with $r>0$ iff
 ${{\partial f}/{\partial r}}(r,\theta,z)<0$.
\item{(d)} The plane field $\xi$ is a positive contact structure
at a point $(0,\theta,z)\in{\cal C}$ iff
the function $h$
as in (a) satisfies the condition $h(0,\theta,z)<0$.
\end{lemma}

The proof of Lemma \ref{BoxChartLemma} is straightforward, and can be found in
\cite{ET:Confoliations}. We will include the slightly
more difficult proof of Lemma
\ref{CylinderChartLemma}.

\begin{proof} {\it ( Lemma \ref{CylinderChartLemma}.)} To prove (a), consider coordinates $(x,y,z)$ in $\cal C$
with $x=r\cos\theta$, $y=r\sin\theta$. In these coordinates we have
$$
d\theta={x\,dy-y\,dx\over x^2+y^2} \hbox{ and hence }
\omega=dz+{f(x,y,z)\over x^2+y^2}(x\,dy-y\,dx).
$$
For a function $h=f/r^2$ we then have $\omega=dz+h(x,y,z)(x\,dy-y\,dx)$,
and therefore $\omega$ is smooth iff the function $h:{\cal C}\to R$
is well defined and smooth.

Having proved (a), the other parts of the lemma follow
by direct calculation.
\end{proof}

The following lemmas
establish the existence of ``purifying operations", which modify a
contamination to eliminate portions of the non-contact locus.  
In particular, we will use these lemmas to remove foliated 
parts of a contamination obtained (in the next section) from
the foliation $\cal F$ by inserting pieces of contact structures
in the complement $V(B)\setminus\pi^{-1}(N(X))$.

\begin{lemma}\label{PurificationLemmaDiskProduct}  Suppose $\xi$ is  
a contamination carried by a branched surface $B\embed M$, and  
let $\psi\from {\cal C}\to V(B)$ be a $\cal C$-chart for $\xi$. 
Identifying $\cal C$ with its image in $V(B)$, suppose
$\xi$ satisfies one of the following two conditions:
\item{(1)} $\xi$ is contact in the subset $\{|z|<1, r<r_0\}\subset{\cal C}$ 
for some $0<r_0<R$;
\item{(2)} $\xi$ is contact in the subset $\{|z|<1, r>r_0\}\subset{\cal C}$
for some $0<r_0<R$.

\noindent
Then $\xi$ can be modified to a contamination $\xi'$ in $V(B)$ (still carried by $B$)
which coincides with $\xi$ outside $\cal C$ and which is pure in the whole
of $\cal C$ except at the top and bottom, i.e. except at $\{z=\pm1\}$.

\end{lemma}

\begin{lemma}\label{PurificationLemmaRectangleProduct} 

Suppose $\xi$ is  a contamination carried by a branched surface $B\embed M$, and
let $\psi\from {\cal B}\to V(B)$ be a $\cal B$-chart for $\xi$. 
Identifying $\cal B$ with its image in $V(B)$, suppose $\xi$ is contact
in the subset $\{|z|<1, |y|>y_0\}$ for some $0<y_0<1$. Then for an arbitrarily small
$\delta>0$, $\xi$ can be modified to a contamination $\xi'$ in $V(B)$ (still carried by $B$)
which coincides with $\xi$ outside $\cal B$ and which is contact in the subset
$\{|z|<1, |x|<1-\delta\}\subset{\cal B}$.

\end{lemma}

We wil prove Lemma 6.7 and omit the proof of Lemma 6.6, which is essentially the same.

\hop

\noindent
{\it Proof of Lemma 6.7.}  Let $f:{\cal B}\to R$ be the slope function for the contamination
$\xi$, in the $\cal B$-chart coordinates as above.
We then have $\partial f/\partial y\le 0$ and $\partial f/\partial y(x,y,z)<0$
for $y>y_0$ and $|z|<1$. Thus if $1>y_1>y_0$ then $f(x,y_1,z)<f(x,-y_1,z)$
for all $x,z\in(-1,1)$. It follows that there exists a smooth function
$f':{\cal B}\to R$ with the following properties:
\item{(1)} $f'$ coincides with $f$ in the union of the subsets 
$\{|x|\ge1-\delta\}$, $\{y\ge y_1\}$ and $\{|z|=1\}$ in $\cal B$;
\item{(2)} $\partial f'/\partial y<0$ in the subset $\{|x|<1-\delta, |z|<1\}$.

\hop

\noindent
Since the function $f'$ coincides with $f$ near the boundary of $\cal B$ in $V(B)$,
it can be used as a slope function of a new contamination $\xi'$ which coincides
with $\xi$ outside $\cal B$. By Lemma 6.5, $\xi'$ is contact in the subset
$\{|x|<1-\delta, |z|<1\}$, which finishes the proof.

\section{Extension and Purification}\label{Purification}

To complete the proof of Theorem 1.1 we will further modify the foliation 
$\cal F$ obtained in Section 5. We will first extend $\cal F$ to an impure
contamination $\xi_0$ in $V(B)$, by inserting standard pieces of contact structure
in the components of $V(B)\setminus\pi^{-1}(N(X))$. Next, we will purify the
contamination $\xi_0$ using the special charts and purifying operations
described in the previous section.
   
Consider the cell structure in $B$ as in Section 4, with its 1-skeleton $X$ 
and with 2-cells denoted by $E$. Recall that $\cal F$ is a smooth foliation in
$\pi^{-1}(N(X))$ transverse to the fibers in $V(B)$. We can choose in each
part $\pi^{-1}(E)$ of $V(B)$ smooth cyllindrical coordinates 
$(r,\theta,z)$, $|z|\le1$, $r\le R$, with the following properties:
\hop
\item{(1)} vertical fibers in $\pi^{-1}(E)$ correspond to curves $(r,\theta)=const$;
\item{(2)} $\pi^{-1}(E\cap N(X))=\{r\ge r_0\}$ for some $0<r_0<R$;
\item{(3)} $\cal F$ is the kernel of a 1-form $dz+f(\theta,z)d\theta$,
where the function $f$ does not depend on $r$.
\hop

\noindent
From the fact that foliation $\cal F$ has strictly negative holonomy
(except at top and bottom), we can further assume (without loss of generality)
that

\item{(4)} $f(\theta,z)<0$ for all $|z|<1$ and all $\theta$.
\hop

It follows from properties of $f$ that there exists a function $f_0$ on $\pi^{-1}(E)$ such that
\item{(a)} $f_0(r,\theta,z)=f(\theta,z)$ for $r\ge r_0$;
\item{(b)} $\partial f_0/\partial r<0$ at $\{|z|<1, 0<r<r_0\}$;
\item{(c)} $f_0=r^2\cdot h$ with $h(0,\theta,z)<0$ for all $|z|<1$.

\noindent
In view of Lemma 6.5, $f_0$ defines an extension of the foliation 
$\cal F$ to the plane field
which is contact in the subset $\{|z|<1, r<r_0\}\subset\pi^{-1}(E)$. Similar extensions of 
$\cal F$ in the other components of 
$V(B)\setminus\pi^{-1}(N(X))$ combine to give a contamination
$\xi_0$ in $V(B)$, carried by $B$, but not pure.

Our goal now is to purify $\xi_0$. We do this
in the following four steps.

{\it Step 1.} Observe that cylindrical charts as above in the parts $\pi^{-1}(E)\subset
V(B)$ are in fact $\cal C$-charts for the contamination $\xi_0$. 
Using these $\cal C$-charts and Lemma 6.6(1), we modify
$\xi_0$ so that, for all $E$, it becomes pure in $\pi^{-1}(E)$.
Thus now the locus of interior points of $V(B)$ where $\xi_0$ is not contact
is contained in the preimage $\pi^{-1}(X)$ of the 1-skeleton $X$.

{\it Step 2.} For each edge $e$ in $B$ not in the branch locus one can consider
a $\cal B$-chart for $\xi_0$ in $V(B)$ 
so that $e$ corresponds to the curve $y=0$ in $(x,y)$ coordinates
of this $\cal B$-chart and so that $\xi_0$ is contact at $\{|z|<1, y\ne0\}\subset {\cal B}$.
Using such charts and Lemma 6.7 we can modify $\xi_0$ so that the locus of interior
points of $V(B)$ where $\xi_0$ is not contact is contained in the  $\pi$-preimage 
of an arbitrarily small neighbourhood in $X$ of the branch locus of $B$.
Some further purification with use of appropriate $\cal C$-charts and Lemma 6.6(1)
makes this locus contained in the preimage of the branch locus of $B$ only.

{\it Step 3.}  
By modifying slope functions in small open
$\cal B$-charts adjacent to the circle components in $\partial_vV(B)$, we can 
modify $\xi_0$ introducing thin annular leaves adjacent to
those circles. This can be done so that the remainder of the locus where $\xi_0$
is contact is unchanged.
We then split $V(B)$ on the above annular leaves 
and obtain a neighborhood $V(B')$ of an
isomorphic branched surface
$B'$ and the induced contamination $\xi_0'$ in $V(B')$. 

{\it Step 4.}  
We can now purify $\xi_0'$ near
the locus where it is not contact inside $V(B')$
(we denote this locus by $L$).
The locus $L$ is contained in the subset of $V(B')$ corresponding (before splitting
$B\to B'$) to the preimage by $\pi$ of the branch locus of $B$. 
Denoting by $\pi':V(B')\to B'$ the natural fiberwise
projection in $V(B')$, observe that $L$ is disjoint 
with the preimage by $\pi'$ of the branch locus of $B'$.

We purify first near the preimage of the regular part of the branch locus in $B$,
using Lemma 6.7 in the same manner as in Step 2.
Then, we purify  near the preimages of the points $P$ in the branch locus of $B$,
where two branching lines meet. This
can be done using Lemma
\ref{PurificationLemmaDiskProduct} with appropriate $\cal C$-charts. 
In this way
we obtain a pure contamination $\xi$ carried by the branched surface $B'$
and, since $B'$ is isomorphic to $B$,
Theorem 1.1 follows.

\bibliographystyle{amsplain}
\bibliography{ReferencesUO2}
\end{document}